%
\def\conv{\mathop{\vrule height2,6pt depth-2,3pt 
    width 5pt\kern-1pt\rightharpoonup}}

\advance\vsize by 1 true cm
%
\def\dess #1 by #2 (#3){
  \vbox to #2{
    \hrule width #1 height 0pt depth 0pt
    \vfill
    \special{picture #3} 
    }
  }

\def\dessin #1 by #2 (#3 scaled #4){{
  \dimen0=#1 \dimen1=#2
  \divide\dimen0 by 1000 \multiply\dimen0 by #4
  \divide\dimen1 by 1000 \multiply\dimen1 by #4
  \dess \dimen0 by \dimen1 (#3 scaled #4)}
  }
%
\def \trait (#1) (#2) (#3){\vrule width #1pt height #2pt depth #3pt}
\def \fin{\hfill
	\trait (0.1) (5) (0)
	\trait (5) (0.1) (0)
	\kern-5pt
	\trait (5) (5) (-4.9)
	\trait (0.1) (5) (0)
\medskip}
%


\font\sevenbf=cmbx7

\baselineskip=15pt
\abovedisplayskip=15pt plus 4pt minus 9pt
\belowdisplayskip=15pt plus 4pt minus 9pt
\abovedisplayshortskip=3pt plus 4pt
\belowdisplayshortskip=9pt plus 4pt minus 4pt
\let\epsilon=\varepsilon

\def\biblio #1 #2\par{\parindent=30pt\item{}\kern -30pt\rlap{[#1]}\kern
30pt #2\smallskip}
 %
\catcode`\@=11
\def\@lign{\tabskip=0pt\everycr={}}
\def\equations#1{\vcenter{\openup1\jot\displ@y\halign{\hfill\hbox
{$\@lign\displaystyle##$}\hfill\crcr
#1\crcr}}}
\catcode`\@=12
%
\def\pmb#1{\setbox0=\hbox{#1}%
\hbox{\kern-.04em\copy0\kern-\wd0
\kern.08em\copy0\kern-\wd0
\kern-.02em\copy0\kern-\wd0
\kern-.02em\copy0\kern-\wd0
\kern-.02em\box0\kern-\wd0
\kern.02em}}
%
\def\undertilde#1{\setbox0=\hbox{$#1$}
\setbox1=\hbox to \wd0{$\hss\mathchar"0365\hss$}\ht1=0pt\dp1=0pt
\lower\dp0\vbox{\copy0\nointerlineskip\hbox{\lower8pt\copy1}}}
%

%

\def\maj#1#2,{\rm #1\sevenrm #2\rm{}}
\def\Maj#1#2,{\bf #1\sevenbf #2\rm{}}
\outer\def\lemme#1#2 #3. #4\par{\medbreak
\noindent\maj{#1}{#2},\ #3.\enspace{\sl#4}\par
\ifdim\lastskip<\medskipamount\removelastskip\penalty55\medskip\fi}

\def\Remark #1. {\noindent{\Maj REMARK,\ \bf #1. }}

\outer\def\Lemme#1#2 #3. #4\par{\medbreak
\noindent\Maj{#1}{#2},\ \bf #3.\rm\enspace{\sl#4}\par
\ifdim\lastskip<\medskipamount\removelastskip\penalty55\medskip\fi}



\def\Notation #1. {\noindent{\Maj NOTATION,\ \bf #1. }}

\def\Example #1. {\noindent{\Maj EXAMPLE,\ \bf #1. }}

\hfuzz=1cm


\catcode`\ˆ=\active     \def ˆ{\`a}
\catcode`\‰=\active     \def ‰{\^a}
\catcode`\=\active     \def {\c c}
\catcode`\Ž=\active    \def Ž{\'e} 

\catcode`\=\active   \def {\`e}
\catcode`\=\active   \def {\^e}
\catcode`\'=\active   \def '{\"e}
\catcode`\"=\active   \def "{\^\i}
\catcode`\•=\active   \def •{\"\i}
\catcode`\™=\active   \def ™{\^o}
\catcode`\š=\active   \defš{}
\catcode`\=\active   \def {\`u}
\catcode`\ž=\active   \def ž{\^u}
\catcode`\Ÿ=\active   \def Ÿ{\"u}
\catcode`\ =\active   \def  {\tau}
\catcode`\¡=\active   \def ¡{\circ}
\catcode`\¢=\active   \def ¢{\Gamma}
\catcode`\¤=\active   \def ¤{\S\kern 2pt}
\catcode`\¥=\active   \def ¥{\puce}
\catcode`\§=\active   \def §{\beta}
\catcode`\¨=\active   \def ¨{\rho}
\catcode`\©=\active   \def ©{\gamma}
\catcode`\­=\active   \def ­{\neq}
\catcode`\°=\active   \def °{\ifmmode\ldots\else\dots\fi}
\catcode`\±=\active   \def ±{\pm}
\catcode`\²=\active   \def ²{\le}
\catcode`\³=\active   \def ³{\ge}
\catcode`\µ=\active   \def µ{\mu}
\catcode`\¶=\active   \def ¶{\delta}
\catcode`\·=\active   \def ·{\Sigma}
\catcode`\¸=\active   \def ¸{\Pi}
\catcode`\¹=\active   \def ¹{\pi}
\catcode`\»=\active   \def »{\Upsilon}
\catcode`\¾=\active   \def ¾{\alpha}
\catcode`\À=\active   \def À{\cdots}
\catcode`\Â=\active   \def Â{\lambda}
\catcode`\Ã=\active   \def Ã{\sqrt}
\catcode`\Ä=\active   \def Ä{\varphi}
\catcode`\Å=\active   \def Å{\xi}
\catcode`\Æ=\active   \def Æ{\Delta}
\catcode`\Ç=\active   \def Ç{\cup}
\catcode`\È=\active   \def È{\cap}
\catcode`\Ï=\active   \def Ï{\oe}
\catcode`\Ñ=\active   \def Ñ{\to}
\catcode`\Ò=\active   \def Ò{\in}
\catcode`\Ô=\active   \def Ô{\subset}
\catcode`\Õ=\active   \def Õ{\superset}
\catcode`\Ö=\active   \def Ö{\over}
\catcode`\×=\active   \def ×{\nu}
\catcode`\Ù=\active   \def Ù{\Psi}
\catcode`\Ú=\active   \def Ú{\Xi}
\catcode`\Ü=\active   \def Ü{\omega}
\catcode`\Ý=\active   \def Ý{\Omega}
\catcode`\ß=\active   \def ß{\equiv}
\catcode`\à=\active   \def à{\chi}
\catcode`\á=\active   \def á{\Phi}
\catcode`\ä=\active   \def ä{\infty}
\catcode`\å=\active   \def å{\zeta}
\catcode`\æ=\active   \def æ{\varepsilon}
\catcode`\è=\active   \def è{\Lambda}  
\catcode`\é=\active   \def é{\kappa}
\catcode`\ë=\active   \defë{\Theta}
\catcode`\ì=\active   \defì{\eta}
\catcode`\í=\active   \defí{\theta}
\catcode`\î=\active   \defî{\times}
\catcode`\ñ=\active   \defñ{\sigma}
\catcode`\ò=\active   \defò{\psi}

\def\date{\number\day\
\ifcase\month \or janvier \or f\'evrier \or mars \or avril \or mai \or juin \or juillet \or ao\^ut  \or
septembre \or octobre \or novembre \or d\'ecembre \fi
\ \number\year}

\font \Ggras=cmb10 at 12pt
\font \ggras=cmb10 at 11pt

\def\sym{\fam\comfam\com}
\font\tensym=msbm10
\font\sevensym=msbm7
\font\fivesym=msbm5
\newfam\symfam
\textfont\symfam=\tensym
\scriptfont\symfam=\sevensym
\scriptfont\symfam=\fivesym
\def\sym{\fam\symfam\relax}
\def\N{{\sym N}}
\def\R{{\sym R}}

\def\Z{{\sym Z}}

\hoffset 0,1cm

\font \Ggras=cmb10 at 14pt
\font \ggras=cmb10 at 16pt

\def\biblio #1 #2\par{\parindent=30pt\item{}\kern -30pt\rlap{[#1]}\kern
30pt #2\smallskip}
\def\biblio #1 #2\par{\parindent=30pt\item{[]}\kern
-30pt\rlap{[#1]}\kern 30pt #2\smallskip}


\centerline {\Ggras Error estimate and unfolding for  periodic homogenization}
\vskip 5mm
\centerline{G. Griso}
\vskip 5mm
\centerline{  Laboratoire J.-L. Lions--CNRS, Bo\^\i te courrier 187, Universit\'e  Pierre et
Marie Curie, }
\centerline{ 4~place Jussieu, 75005 Paris, France,  Email: griso@ann.jussieu.fr}
\vskip 8mm
\noindent{\ggras \bf Abstract. }
{\sevenrm This paper deals with the error estimate in problems of periodic homogenization. The methods
used  are those of the periodic unfolding. We give the upper bound  of the distance between the unfolded gradient of a 
function belonging to $\scriptstyle H^1(\Omega)$ and the space  $\scriptstyle \nabla_x H^1(\Omega)\oplus \nabla_y
L^2(\Omega ; H^1_{per}(Y))$. These distances  are obtained thanks to a technical result presented in Theorem 2.3 : the periodic
defect of a harmonic function belonging to $\scriptstyle H^1(Y)$ is written with the help of the norms $\scriptstyle H^{1/2}$ of
its  traces differences on the opposite faces of the cell $\scriptstyle Y$. The error estimate is obtained without any
supplementary hypothesis of regularity on correctors. }
\bigskip
\noindent {\bf 1. Introduction}
\vskip 1mm
\noindent The error estimate in periodic homogenization problems was presented for the first time in  Bensoussan, Lions
and Papanicolaou [2]. It can also be found     in   Oleinik, Shamaev and Yosifian  [8], and more recently   in    Cioranescu
and Donato [5].  In all these books, the result  is proved under the assumption that the correctors belong to $W^{1,\infty}(Y)$
($Y=]0,1[^n$ being the reference cell). The  estimate is of order   $\varepsilon^{1/ 2}$. The  additional  regularity of the
correctors holds true  when the coefficients of the operator are very regular, which is not necessarily the situation in
homogenization. In [6]  we  obtained an error estimate without any regularity hypothesis on the correctors but we supposed
that the solution of the homogenized problem belonged to
$W^{2,p}(\Omega)$ ($p>n$). The exponent of $\varepsilon$ in the error estimate is inferior to $1/2$ and depends on $n$ and $p$. 

The aim of this work is to give further error estimates with again minimal hypotheses on the correctors and the homogenized
problem. In all this  study we will make use of  the notation of [4].
 
The paper is organized as follows. In paragraph 2  we prove some technical results related to periodic defect.  In Theorem 2.1 we
give an estimate of the distance between a function $\phi$ belonging to $W^{1,p}(Y)$ and the space of periodic functions
$W^{1,p}_{per}(Y)$. This distance depends on the $W^{1-{1\over p},p}$ norms of the differences of the traces of $\phi$ on
opposite faces of $Y$. Theorem 2.1 is  a consequence of Lemma 2.2. In this lemma we proved that the distance between a
function  and the space of   periodic functions with respect to the  first $k$ variables   is isomorphic to the  direct sum of the
spaces of the differences of the traces on the opposite faces  $Y_j$ and $\vec e_j+Y_j$.  (i.e., on 
$y_{j}=0\;\hbox{and}\; y_{j}=1$), $1\le j\le k$.  This lemma is proved by an explicit lifting of the traces from the faces of $Y$.  

In Theorem  2.3 we show that the $H^{1/2}$ periodic defect of an harmonic function on $Y$ with values in a separable Hilbert
space $X$ is equivalent to its $H^1$ norm. The orthogonal of  space $H^1_{per}(Y ; X)$ is in fact isomorphic to the direct sum of
the spaces of the differences of the traces on the opposite faces of  cell $Y$. 

Paragraph 3 is dedicated to  Theorem 3.4 which is the essential tool  to obtain estimates. This theorem is related to  the periodic
unfolding method (see [4]).   We show that for any $\phi$ in $H^1(\Omega)$, where $\Omega$ is an open bounded set of
$\R^n$  with Lipschitz boundary, there exists a function $\widehat{\phi}_\varepsilon$ in $H^1_{per}(Y ; L^2(\Omega))$, such
that the distance between the unfolded ${\cal T}_\varepsilon(\nabla_x\phi)$
and $\nabla_x\phi+\nabla_y\widehat{\phi}_\varepsilon$ is of order of
$\varepsilon$ in the space $[L^2(Y ; H^{-1}(\Omega))]^n$. 

 Theorems 4.1, 4.2 and 4.5 give an estimate of the  error without any hypothesis on the regularity of the correctors, but with
different hypotheses on the boundary of $\Omega$. They require that the right hand side of the homogenized problem be in  
$L^2(\Omega)$.

In this article, the constants appearing in the estimates will be independent from $\varepsilon$.
\vfill\eject
\noindent{\bf 2. Periodicity defect}
\medskip
\noindent  We denote  $Y=]0,1[^n$ the unit cell of $\R^n$ and we put $Y_j=\bigl\{y\in \overline{Y}\; |\; y_j=0\bigr\}$,
$j\in\{1,\ldots, n\}$.
\medskip
\noindent{\bf Theorem 2.1 : }{\it For any $\phi\in W^{1,p}(Y )$,
$p\in ]1,\infty]$, there exists  $\widehat{\phi}\in W^{1,p}_{per}(Y )$ such that
$$||\phi-\widehat{\phi}||_{W^{1,p}(Y)}\le C\sum_{j=1}^n||\phi_{|_{\vec e_j+Y_j}}-\phi_{|_{Y_j}}||_{W^{1-{1\over p},p}(Y_j
)}$$ The constant depends only on $n$.}

\noindent  The proof of the theorem is based on Lemma 2.2.  We introduce the following spaces:
$$W_0=W^{1,p}(Y),\qquad\quad W_k  =\Bigl\{\phi\in W^{1,p}(Y )\; |\; \phi(.)=\phi(.+\vec e_i),\enskip
i\in\{1,\ldots,k\}\Bigr\},\quad k\in\{1,\ldots,n\}$$
\noindent{\bf Lemma 2.2 : }{\it For any $\phi\in W^{1,p}(Y)$ and for any $k\in \{1,\ldots, n\}$,  there exists
$\widehat{\phi}_k\in W_k$ such that
$$ ||\phi-\widehat{\phi}_k||_{W^{1,p}(Y)}\le C \sum_{j=1}^k||\phi_{|_{\vec e_j+Y_j}}-\phi_{|_{Y_j}}||_{W^{1-{1\over
p},p}(Y_i )}$$ The constant depends  on $n$.}
\medskip
\noindent{\bf Proof : } The lemma is proved by  a finite induction. We choose a function $\theta$ belonging to
${\cal D}(-1/2,1/2)$ equal to 1 in the neighborhood of zero.  We recall (see [1]) that for any  $k\in\{0,\ldots, n-1\}$, there
exists a continuous lifting $\widetilde{ r}_k$ in $W_k$ of the traces on $Y_{k+1}$ of the $W_k$ elements. 

\noindent Let $\phi$ be in $W^{1,p}(Y)$. We put $\widehat{\phi}_0=\phi$. We suppose  the lemma proved for $k$,
$k\in\{0,\ldots, n-1\}$. There exists $\widehat{\phi}_k\in W_k$ such that
$$ ||\phi-\widehat{\phi}_k||_{W^{1,p}(Y)}\le C \sum_{j=1}^k||\phi_{|_{\vec e_j+Y_j}}-\phi_{|_{Y_j}}||_{W^{1-{1\over
p},p}(Y_i )}$$ Of course if $k=0$ the right hand side of the above inequality is equal to zero. We define $\widehat{\phi}_{k+1}$
by $$\widehat{\phi}_{k+1}=\widehat{\phi}_k+{1\over
2}\Bigl\{\theta(y_{k+1})-\theta(1-y_{k+1})\Bigr\}\widetilde{r}_k\bigl(\widehat{\phi}_{k|_{\vec
e_{k+1}+Y_{k+1}}}-\widehat{\phi}_{k|_{Y_{k+1}}}\bigr) $$ The function
$\widehat{\phi}_{k+1}$ belongs to $W_k$ and verifies 
$$\widehat{\phi}_{k+1|_{\vec e_{k+1}+Y_{k+1}}}=\displaystyle {1\over 2}\bigl\{\widehat{\phi}_{k|_{\vec
e_{k+1}+Y_{k+1}}}+\widehat{\phi}_{k|_{Y_{k+1}}}\bigr\}=\widehat{\phi}_{k+1|_{Y_{k+1}}}$$
Hence it belongs to $W_{k+1}$. We have
$$ \eqalign{
||\phi-\widehat{\phi}_{k+1}||_{W^{1,p}(Y)} &\le  ||\phi-\widehat{\phi}_k||_{W^{1,p}(Y)}+
||\widehat{\phi}_k-\widehat{\phi}_{k+1}||_{W^{1,p}(Y)}\cr
&\le  C\sum_{j=1}^k||\phi_{|_{\vec e_j+Y_j}}-\phi_{|_{Y_j}}||_{W^{1-{1\over p},p}(Y_i )}+ C||\widehat{\phi}_{k|_{\vec
e_{k+1}+Y_{k+1}}}-\widehat{\phi}_{k|_{Y_{k+1}}}||_{W^{1-{1\over p},p}(Y_{k+1} )}\cr}$$ Besides, we have
$$\eqalign{
||\widehat{\phi}_{k|_{\vec e_{k+1}+Y_{k+1}}}-\widehat{\phi}_{k|_{Y_{k+1}}}||_{W^{1-{1\over p},p}(Y_{k+1} )}&\le 
||(\phi-\widehat{\phi}_k)_{|_{\vec e_{k+1}+Y_{k+1}}}-(\phi-\widehat{\phi}_k)_{|_{Y_{k+1}}}||_{W^{1-{1\over p},p}(Y_{k+1}
)}\cr
&+||\phi_{|_{\vec e_{k+1}+Y_{k+1}}}-\phi_{|_{Y_{k+1}}}||_{W^{1-{1\over p},p}(Y_{k+1})}\cr
&\le C||\phi-\widehat{\phi}_{k}||_{W^{1,p}(Y)}+||\phi_{|_{\vec
e_{k+1}+Y_{k+1}}}-\phi_{|_{Y_{k+1}}}||_{W^{1-{1\over p},p}(Y_{k+1})}\cr}$$ Hence we obtain the result for $k+1$ and the
lemma is proved.\fin
\medskip
\noindent{\bf Proof of Theorem 2.1 : } We have $W^n=W^{1,p}_{per}(Y)$. Thanks to Lemma 2.2   Theorem 2.1
is proved by taking   $k=n$.\fin 
\noindent{\bf Remark 1 : } Let $X$ be a Banach space.  We can prove as in Lemma 2.2 and Theorem 2.1  that
for any $\Phi\in W^{1,p}(Y; X)$ there exists $\widehat{\Phi}\in W^{1,p}_{per}(Y ; X)$ such that
$$||\Phi-\widehat{\Phi}||_{W^{1, p}(Y; X)}\le C\sum_{j=1}^n||\Phi_{|_{\vec e_j+Y_j}}-\Phi_{|_{Y_j}}||_{W^{1-{1\over
p},p}(Y_j; X)}$$ The constant depends only on $n$. \fin
\medskip
\noindent Let $X$ be a separable Hilbert space.  We equip $H^1(Y; X)$  with the inner product
$$<\phi,\psi>=\int_Y \nabla\phi\cdot \nabla\psi +\Bigl(\int_Y\phi\Bigr)\cdot\Bigl(\int_Y\psi\Bigr)  $$ where 
$\cdot$ is the inner product in $X$. The norm associated to this scalar product is equivalent to the norm of $H^1(Y; X)$.

\noindent{\bf Theorem 2.3 : }{\it For any $\phi\in H^1(Y ; X)$ there exists  a unique $\widehat{\phi}\in H^1_{per}(Y ; X)$ such
that
$$\phi-\widehat{\phi}\in \bigl(H^1_{per}(Y; X)\bigr)^\perp\qquad
||\widehat{\phi}||_{H^1(Y ; X)}\le||\phi||_{H^1(Y ; X)}\qquad ||\phi-\widehat{\phi}||_{H^1(Y ; X)}\le
C\sum_{j=1}^n||\phi_{|_{\vec e_j+Y_j}}-\phi_{|_{Y_j}}||_{H^{1/2}(Y_j ; X)}$$ The constant depends only on $n$. The function
$\widehat{\phi}$ verifies
$$\int_Y\phi=\int_Y\widehat{\phi},\qquad\qquad \forall\psi\in \bigl(H^1_{per}(Y ; \R)\bigr)^\perp\qquad
\int_Y\nabla(\phi-\widehat{\phi})\nabla\psi=0\qquad\hbox{in}\quad X.$$}
\noindent{\bf Proof : } We take a Hilbert basis $\bigl(x_n\bigr)_{n\in \N}$ of  $X$. Any element $\phi$ belonging to $H^1(Y;
X)$ is decomposed into a series $\displaystyle \phi=\sum_{n=0}^\infty \phi_n x_n$, where $\phi_n$ belongs to $H^1(Y)$. 
We apply  Theorem 2.1 to each component $\phi_n$ and then by orthogonal projection we obtain Theorem 2.3.\fin

\noindent{\bf Corollary : } If ${\cal X}$ is a Hilbert space continuously embedded in $X$ then for any $\phi \in H^1(Y;{\cal
X})$ there exists $ \widehat{\phi}\in H^1_{per}(Y;{\cal X} )$ such that
$$\eqalign{
||\phi-\widehat{\phi}||_{H^1(Y;{\cal X})}\le C\sum_{j=1}^n||\phi_{|_{\vec e_j+Y_j}}-\phi_{|_{Y_j}}||_{H^{1/2}(Y_j;{\cal X} )}\cr
||\phi-\widehat{\phi}||_{H^1(Y;X)}\le C\sum_{j=1}^n||\phi_{|_{\vec e_j+Y_j}}-\phi_{|_{Y_j}}||_{H^{1/2}(Y_j ; X )}\cr}$$ 
The constant depends only on $n$.
\medskip
\noindent{\bf 3. Approximation and periodic unfolding}
\vskip 1mm
\noindent Let  $\Omega$   be a bounded domain in $\R^n$ with lipschitzian boundary. We put
$$\eqalign{
&\widetilde{\Omega}_{\varepsilon, k}=\bigl\{x\in \R^n\; | \; dist(x,  \Omega)< k\sqrt n\, \varepsilon\;\bigr\},\qquad
k\in\{1,2\},\cr
&\Omega_\varepsilon =\hbox{Int}\Bigl(\bigcup_{\xi\in \Xi_\varepsilon}\varepsilon(\xi+\overline{Y})\Bigr),\qquad
\Xi_\varepsilon=\Bigl\{\xi\in \Z^n\; |\; \varepsilon(\xi+\overline{Y})\cap\Omega\not=\emptyset\Bigr\}\cr}$$  We have 
$$\Omega\i \Omega_\varepsilon\i \widetilde{\Omega}_{\varepsilon,1},\qquad \hbox{and}\qquad\forall i\in\{1,\ldots,
n\},\qquad \Omega_\varepsilon+\varepsilon\vec e_i\i \widetilde{\Omega}_{\varepsilon,2}.$$
\noindent We recall that there exists a linear and continuous extension operator ${\cal P}$ from
$H^1(\Omega)$  into $H^1(\widetilde{\Omega}_{\varepsilon, 2})$, such that for any
$\phi\in H^1(\Omega)$, ${\cal P}(\phi)$ belongs to $H^1(\widetilde{\Omega}_{\varepsilon,2})$ and verifies 
$$\eqalign{
 &{\cal P}(\phi)_{|_\Omega}=\phi,\qquad ||\nabla_x {\cal P}(\phi)||_{[L^2(\widetilde{\Omega}_{\varepsilon,2})]^n}\le
C||\nabla_x\phi||_{[L^2(\Omega)]^n}\cr  &||{\cal P}(\phi)||_{L^2(\widetilde{\Omega}_{\varepsilon,2})}\le C\bigl\{||\phi
||_{L^2(\Omega)}+\varepsilon||\nabla_x\phi||_{[L^2(\Omega)]^n}\bigr\}\cr}$$ More precisely, we have
$$||{\cal P}(\phi)||_{L^2(\widetilde{\Omega}_{\varepsilon,2})}+\varepsilon||\nabla_x {\cal
P}(\phi)||_{[L^2(\widetilde{\Omega}_{\varepsilon,2})]^n}\le C\bigl\{||\phi
||_{L^2(\Omega)}+\varepsilon||\nabla_x\phi||_{[L^2(\Omega)]^n}\bigr\}\leqno(3.1)$$  {\sl In the rest of this paragraph,
 without having to specify it every time, any function belonging to $H^1(\Omega)$ is extended to 
$\widetilde{\Omega}_{\varepsilon,2}$, the extension verifying $(3.1)$. In order to simplify the notation, we will still denote 
by $\phi$ its extension. }
\medskip
\noindent In the sequel, we will make use of definitions and results from [4] concerning the periodic unfolding method.  

\noindent  For almost every $x$ belonging to $ \R^n$, there exists a unique
element in $\Z^n$ denoted  $[x]$  such that 
$$x=[x]+\{x\},\qquad \{x\}\in Y.$$ 
Let us  now recall the definition  of the unfolding operator ${\cal T}_\varepsilon$ which to each function  $\phi\in  L^1(
\Omega_\varepsilon)$   associates a function ${\cal T}_\varepsilon(\phi)\in L^1(\Omega\times Y)$, 
$${\cal T}_\varepsilon(\phi)(x,y)=\phi\Bigl(\varepsilon\Bigr[{x\over \varepsilon}\Bigr]+\varepsilon y\Bigr)\qquad \hbox{ for
$x\in \Omega$ and  $y\in Y$}.$$ We have
$$\Bigl|\int_\Omega \phi-\int_{\Omega\times Y}{\cal T}_\varepsilon(\phi)\Bigr|\le
||\phi||_{L^1( \{x\in \Omega_\varepsilon\; | \; dist(x,  \partial\Omega)< \sqrt n \,\varepsilon\; \})}$$
\noindent For the other properties of  ${\cal T}_\varepsilon$, we refer the reader to  [4]. 
\medskip
\noindent  Now, for any  $\phi\in L^2(\Omega_\varepsilon)$ we define the operator  ``mean in the cells''
$M^\varepsilon_Y$ by setting 
$$\displaystyle M^\varepsilon_Y(\phi)(x)=\int_Y{\cal T}_\varepsilon(\phi)(x,y)dy={1\over \varepsilon^n}\int_{\displaystyle
\bigl\{\varepsilon\bigl[{x\over\varepsilon}\bigr]+\varepsilon Y\bigr\}}\phi(z)dz,\qquad x\in\Omega.$$ Function 
$M^\varepsilon_Y(\phi)$ belongs to $L^2(\Omega)$ and verifies
$$||M^\varepsilon_Y(\phi)||_{L^2(\Omega)}\le ||\phi||_{L^2(\Omega_\varepsilon)}.$$
\noindent{\bf Proposition 3.1 : }{\it  For any $\phi$ belonging to $H^1(\Omega)$ we have
$$||\phi-M^\varepsilon_Y(\phi)||_{L^2(\Omega)}\le C\varepsilon||\nabla_x\phi||_{[L^2(\Omega)]^n}.\leqno(3.2)$$}
\noindent{\bf Proof : } Let $\phi\in H^1(\Omega)$. We apply the Poincar\'e-Wirtinger inequality to the  restrictions
$x\longrightarrow \phi_{|_{\varepsilon(\xi+Y)}}(x)-M^\varepsilon_Y(\phi)(\varepsilon\xi)$ belonging to
$H^1(\varepsilon(\xi+Y))$
$$||\phi-M^\varepsilon_Y(\phi)(\varepsilon\xi)||^2_{L^2(\varepsilon(\xi+Y))}\le
C\varepsilon^2||\nabla_x\phi||^2_{[L^2(\varepsilon(\xi+Y))]^n},\qquad
\varepsilon(\xi+Y)\i \Omega_\varepsilon$$ We add all these inequalities and obtain  $(3.2)$.\fin

\noindent We recall the definition of the scale-splitting operator ${\cal Q}_\varepsilon$. The function ${\cal Q}_\varepsilon(\phi)$
is  the  restriction to $\Omega$ of  $Q_1$-interpolate of the discrete function $M^\varepsilon_Y(\phi)$.

\noindent{\bf Corollary : }{\it For any $\phi\in L^2(\Omega_\varepsilon)$ we have
$$||\phi-M^\varepsilon_Y(\phi)||_{H^{-1}(\Omega)}\le C\varepsilon||\phi||_{L^2(\Omega_\varepsilon)}.\leqno(3.3)$$
\noindent For any $\phi\in H^1(\Omega)$ we have
$$\left\{\eqalign{  & ||\phi-{\cal T}_\varepsilon(\phi)||_{L^2(\Omega\times Y)}\le
C\varepsilon||\nabla_x\phi||_{[L^2(\Omega)]^n}.\cr & ||{\cal Q}_\varepsilon(\phi)-M^\varepsilon_Y(\phi)||_{L^2(\Omega)}\le
C\varepsilon||\nabla_x\phi||_{[L^2(\Omega)]^n}.\cr}\right.\leqno(3.4)$$}
\noindent{\bf Proof : } If $\psi\in H^1_0(\Omega)$, we immediately have
$$\int_\Omega\bigl(\phi-M^\varepsilon_Y(\phi)\bigr)\psi=\int_{\Omega_\varepsilon}\bigl(\phi-M^\varepsilon_Y(\phi)\bigr)\psi
=\int_{\Omega_\varepsilon}\phi\bigl(\psi- M^\varepsilon_Y(\psi)\bigr)\le C\varepsilon
||\phi||_{L^2(\Omega_\varepsilon)}||\nabla_x\psi||_{[L^2(\Omega)]^n}$$ hence   inequality  $(3.3)$.

\noindent We have (see [4]): if $\phi\in L^2(\Omega_\varepsilon)$ then
$$||{\cal T}_\varepsilon \bigl(\phi-M^\varepsilon_Y(\phi)\bigr) ||_{L^2(\Omega\times Y)}
\le ||\phi-M^\varepsilon_Y(\phi)||_{L^2(\Omega_\varepsilon)}$$ and moreover,  ${\cal T}_\varepsilon\circ
M^\varepsilon_Y(\phi)=M_Y^\varepsilon(\phi)$.  We eliminate the mean function  $M^\varepsilon_Y(\phi)$ with $(3.2)$  to
obtain $(3.4)$.  We also have (see [4]) $ ||\phi-{\cal Q}_\varepsilon(\phi)||_{L^2(\Omega)}\le
C\varepsilon||\nabla_x\phi||_{[L^2(\Omega)]^n}$ and according to 
$(3.2)$ we obtain the second inequality of $(3.4)$. \fin
\noindent{\bf Proposition 3.2 : }{\it For any  $\phi$ belonging to $ L^2(\widetilde{\Omega}_{\varepsilon, 2})$ and any $\psi$
belonging to $L^2(Y)$, we have
$$||{\cal Q}_\varepsilon(\phi)\psi\bigl(\bigl\{{.\over \varepsilon}\bigr\}\bigr)||_{L^2(\Omega)}\le
C||\phi||_{L^2(\widetilde{\Omega}_{\varepsilon, 2})}||\psi||_{L^2(Y)}\leqno(3.5)$$ The constant depends only on $n$.}

\noindent{\bf Proof : }  We set for 
$i=(i_1,\ldots, i_n)\in \{0,1\}^n$,
$$x\in \varepsilon\bigl(\xi+Y),\qquad\qquad \overline{x}_{i, \xi}^{i_k}=\left\{\eqalign{ 
&{x_k-\varepsilon\xi_k\over \varepsilon}\qquad\quad\;\;\hbox{if}\enskip i_k=1
\cr &1- {x_k-\varepsilon\xi_k\over \varepsilon}\qquad\hbox{if}\enskip i_k=0 \cr}\right.$$ From the definition of
${\cal Q}_\varepsilon(\phi)$ (see [4]) it results  that
$$x\in \Omega_\varepsilon,\qquad\qquad {\cal Q}_\varepsilon(\phi)(x)=\sum_{i_1\, \ldots, i_n}M^\varepsilon_Y(\phi)\bigl(
\varepsilon\xi+\varepsilon i\bigr)\overline{x}_{1 , \xi}^{i_1}\ldots \overline{x}_{ n , \xi}^{i_n},\qquad \xi=\Bigl[{x\over
\varepsilon}\Bigr]$$ hence
$$\eqalign{
\int_{\varepsilon(\xi+Y)} |{\cal Q}_\varepsilon(\phi)|^2|\psi\bigl(\bigl\{{.\over \varepsilon}\bigr\}\bigr)|^2 &\le 2^n\sum_{i_1\, \ldots, i_n}
|M^\varepsilon_Y(\phi)\bigl(\varepsilon\xi+\varepsilon i\bigr)|^2\int_{\varepsilon(\xi+Y)}|\psi\bigl(\bigl\{{.\over \varepsilon}\bigr\}\bigr)|^2\cr
 &= 2^n\sum_{i_1\, \ldots, i_n} |M^\varepsilon_Y(\phi)\bigl(\varepsilon\xi+\varepsilon
i\bigr)|^2\varepsilon^n ||\psi||^2_{L^2(Y)}\cr}$$  For any $\xi $ we have  $|M^\varepsilon_Y(\phi)(\varepsilon\xi)|^2\le
\displaystyle{1\over \varepsilon^n|Y|}\int_{\varepsilon (\xi+Y)}|\phi|^2$. We add the above inequalities  for all $\xi\in
\Xi_\varepsilon$ and we obtain
$$\int_{\Omega} |{\cal Q}_\varepsilon(\phi)|^2|\psi\bigl(\bigl\{{.\over \varepsilon}\bigr\}\bigr)|^2\le 4^n
||\phi||^2_{L^2(\widetilde{\Omega}_{\varepsilon, 2})}||\psi||^2_{L^2(Y)}$$\fin
\noindent{\bf Proposition 3.3 : }{\it For any  $\phi$ belonging to $H^1(\Omega)$, there exists  $\widehat{\psi}_\varepsilon$
belonging to $H^1_{per}(Y ; L^2(\Omega))$  \footnote{$^{(1)}$ }{Of course $H^1_{per}(Y; L^2(\Omega))$ is the same space
as $L^2(\Omega; H^1_{per}(Y))$. The same remark holds for all other spaces appearing in the sequel.} such that
$$\left\{\eqalign{ &||\widehat{\psi}_\varepsilon||_{H^1(Y ; L^2(\Omega))}\le
C\bigl\{||\phi||_{L^2(\Omega)}+\varepsilon||\nabla_x\phi||_{[L^2(\Omega)]^n}\bigr\}\cr &||{\cal
T}_\varepsilon(\phi)-\widehat{\psi}_\varepsilon||_{ H^1(Y ; H^{-1}(\Omega ))}\le C\varepsilon
\bigl\{||\phi||_{L^2(\Omega)}+\varepsilon||\nabla_x\phi||_{[L^2(\Omega)]^n}\bigr\}\cr}\right.\leqno(3.6)$$}
\noindent{\bf Proof : } Proposition 3.3 is proved in two steps. We begin with constructing a new unfolding operator which for any
$\phi\in H^1(\Omega)$ allows us to estimate in
$L^2(Y ; H^{-1}(\Omega))$, the difference between the restrictions to two neighboring cells of the unfolded of
$\phi$. Then, we  evaluate the periodic defect of the functions   $y\longrightarrow {\cal T}_\varepsilon(\phi)(.,y)$ and
conclude thanks  to Theorem 2.3.
\vskip 1mm
\noindent Let  $\displaystyle K_i=\hbox{Int}\bigl(\overline{Y}\cup (\vec e_i+\overline{Y})\bigr)$, $i\in\{1,\ldots,n\}$. For any
$x$ in
$\Omega$,
$\displaystyle 
\varepsilon\Bigl(\Bigl[{x\over \varepsilon}\Bigr]+ K_i\Bigr)$ is included in $\widetilde{\Omega}_{\varepsilon,2}$.
\vskip 1mm
\noindent{\bf Step one.}  We define the unfolding operator
${\cal T}_{\varepsilon, i}$ from $L^2(\widetilde{\Omega}_{\varepsilon,2})$ into $L^2(\Omega\times K_i)$ by
$$\forall\psi\in L^2(\widetilde{\Omega}_{\varepsilon,2}),\qquad {\cal
T}_{\varepsilon,i}(\psi)(x,y)=\psi\Bigl(\varepsilon\Bigl[{x\over
\varepsilon}\Bigr]+ \varepsilon y\Bigr)\qquad \hbox{for $ x\in\Omega$ and  a. e. $ y\in  K_i$}.$$ 
\noindent  The restriction of  ${\cal T}_{
\varepsilon,i}(\psi)$ to $\Omega\times Y$ is equal to the unfolded  ${\cal T}_\varepsilon(\psi)$. Moreover, we have the
following equalities in $L^2(\Omega\times Y)$:
$${\cal T}_{\varepsilon,i}(\psi)(.,..+\vec e_i)={\cal T}_{\varepsilon}(\psi)(.+\varepsilon \vec e_i,..),\qquad
i\in\{1,\ldots,n\}$$ Let us take $\Psi\in H^1_0(\Omega)$, extended by 0 on $\R^n\setminus\Omega$. A linear change of
variables and the above relations give
$$\eqalign{
\hbox{for a. e. } y\in Y,\qquad \int_{\Omega }{\cal T}_{\varepsilon,i}(\psi)(x,y+\vec e_i)\Psi(x)dx&=\int_{\Omega }{\cal
T}_{\varepsilon,i}(\psi)(x+\varepsilon\vec e_i,y)\Psi(x)dx\cr &= \int_{ \Omega+\varepsilon\vec e_i}{\cal
T}_{\varepsilon,i}(\psi)(x,y)\Psi(x-\varepsilon\vec e_i)dx\cr}$$ We deduce 
$$\eqalign{  &\Bigl|\int_{\Omega }\bigl\{{\cal T}_{\varepsilon,i}(\psi)(., y+\vec e_i)-{\cal T}_{\varepsilon,i}(\psi)(.,
y)\bigr\}\Psi
 -\int_{\Omega }{\cal T}_{\varepsilon,i}(\psi)(., y)\bigl\{\Psi(.-\varepsilon\vec e_i)-\Psi\bigr\}\Bigr|\cr
\le & C||{\cal T}_{\varepsilon,i}(\psi)(.,
y)||_{L^2(\widetilde{\Omega}_{\varepsilon,1})}||\Psi||_{L^2(\Omega\Delta\{\Omega+\varepsilon\vec e_i\})}\cr}$$ where
$\Omega\Delta\{\Omega+\varepsilon\vec e_i\}=(\Omega\setminus\{\Omega+\varepsilon\vec e_i\})
\cup (\{\Omega+\varepsilon\vec e_i\}\setminus \Omega)$; $\Omega$ is a bounded domain with lipschitzian boundary and
$\Psi$ belongs to $ H^1_0(\Omega)$, we thus have
$$\eqalign{  &||\Psi||_{L^2(\Omega\Delta\{\Omega+\varepsilon\vec e_i\} )}\le C\varepsilon ||\nabla_x\Psi||_{[L^2(\Omega
)]^n},\cr  & ||\Psi(.-\varepsilon\vec e_i)-\Psi||_{L^2(\Omega)}\le C\varepsilon\Bigl\|{\partial\Psi\over \partial
x_i}\Bigr\|_{L^2(\Omega)},\qquad i\in\{1,\ldots,n\},\cr}$$ hence
$$\eqalign{ &<{\cal T}_{\varepsilon,i}(\psi)(., y+\vec e_i)-{\cal T}_{\varepsilon,i}(\psi)(., y)\,,\,
\Psi>_{H^{-1}(\Omega),H^1_0(
\Omega)}\cr =&\int_{\Omega }\bigl\{{\cal T}_{\varepsilon,i}(\psi)(., y+\vec e_i)-{\cal T}_{\varepsilon,i}(\psi)(.,
y)\bigr\}\Psi\cr
\le  &C\varepsilon||\nabla_x\Psi||_{[L^2(\Omega )]^n}||{\cal T}_{\varepsilon,i}(\psi)(.,
y)||_{L^2(\widetilde{\Omega}_{\varepsilon,1})}
\le  C\varepsilon||\Psi||_{H^1_0(\Omega))}||{\cal T}_{\varepsilon,i}(\psi)(.,
y)||_{L^2(\widetilde{\Omega}_{\varepsilon,1})}.\cr}$$ We deduce
 that
$$||{\cal T}_{\varepsilon,i}(\psi)(.,y +\vec e_i)-{\cal T}_{\varepsilon,i}(\psi)(., y)||_{H^{-1}(\Omega)}\le  C
\varepsilon||{\cal T}_{\varepsilon,i}(\psi)(., y)||_{L^2(\widetilde{\Omega}_{\varepsilon,1}),}$$
\noindent which leads to the following estimate of the difference between ${\cal T}_{\varepsilon,i}(\psi)_{|_{\Omega\times
Y}}$ and one of its translated:
$$||{\cal T}_{\varepsilon,i}(\psi)(.,.. +\vec e_i)-{\cal T}_{\varepsilon,i}(\psi)||_{L^2(Y ; H^{-1}(\Omega))}
\le C\varepsilon||\psi||_{L^2(\widetilde{\Omega}_{\varepsilon,2})}\leqno(3.7)$$ The constant depends only on the boundary of 
$\Omega$.

\noindent{\bf Step two. } Let $\phi\in H^1(\Omega)$. The estimate $(3.7)$ applied to $\phi$ and its partial derivatives gives 
$$\eqalign{ 
||{\cal T}_{\varepsilon,i}(\phi)(., .. +\vec e_i)-{\cal T}_{\varepsilon,i}(\phi) ||_{ L^2(Y ; H^{-1} (\Omega))} & \le
C\varepsilon\bigl\{||\phi||_{L^2(\Omega)}+\varepsilon ||\nabla_x\phi||_{[L^2(\Omega)]^n}\bigr\}\cr  
||{\cal T}_{\varepsilon,i}(\nabla_x\phi)(., .. +\vec e_i)-{\cal T}_{\varepsilon,i}(\nabla_x \phi)||_{ [L^2(Y ; H^{-1}  (\Omega)]^n)} &\le C\varepsilon||\nabla_x\phi||_{[L^2(\Omega)]^n}\cr}$$ We recall
(see [4]) that  $\nabla_y\bigl({\cal T}_{\varepsilon,i}(\phi)\bigr)=\varepsilon{\cal T}_{\varepsilon,i}(\nabla_x\phi)$. The
above estimates can also be written:
$$||{\cal T}_{\varepsilon,i}(\phi)(., .. +\vec e_i)-{\cal T}_{\varepsilon,i}(\phi)||_{H^1(Y ; H^{-1}(\Omega))}\le
C\varepsilon\bigl\{||\phi||_{L^2(\Omega)}+\varepsilon ||\nabla_x\phi||_{[L^2(\Omega)]^n}\bigr\}$$ 
\noindent From these inequalities, for any
$i\in\{1,\ldots,n\}$,  we deduce the estimate of the difference of the traces of
 $y\longrightarrow {\cal T}_{\varepsilon}(\phi)(.,y)$ on the faces $Y_i$ and $\vec e_i+Y_i$ 
$$||{\cal T}_{\varepsilon}(\phi)(.,.. +\vec e_i)-{\cal T}_{\varepsilon}(\phi)||_{H^{1/ 2}(Y_i  ; H^{-1}(\Omega))}\le
C\varepsilon\bigl\{||\phi||_{L^2(\Omega)}+\varepsilon ||\nabla_x\phi||_{[L^2(\Omega)]^n}\bigr\}\leqno(3.8)$$ which measures
the periodic defect of
$y\longrightarrow{\cal T}_\varepsilon(\phi)(.,y)$. Thanks to Theorem 2.3 we decompose  ${\cal T}_{\varepsilon}(\phi)$ in
the sum of an element $\widehat{\psi}_\varepsilon$ belonging to $H^1_{per}(Y ; L^2(\Omega ))$ and an element
$\overline{\phi}_\varepsilon$ belonging to $\bigl(H^1(Y ; L^2(\Omega))\bigr)^\perp$ such that
$$\left\{\eqalign{
||\overline{\phi}_\varepsilon||_{H^1(Y ; H^{-1}(\Omega))} &\le C\sum_{j=1}^n||{\cal T}_{\varepsilon}(\phi)(.,.. +\vec e_j)-{\cal
T}_{\varepsilon}(\phi)||_{H^{1/ 2}(Y_j ; H^{-1}(\Omega))}\cr
&\le C\varepsilon\bigl\{||\phi||_{L^2(\Omega)}+\varepsilon ||\nabla_x\phi||_{[L^2(\Omega)]^n}\bigr\}\cr
||\overline{\phi}_\varepsilon||_{H^1(Y ; L^2(\Omega))} &\le 
C\bigl\{||\phi||_{L^2(\Omega)}+\varepsilon||\nabla_x\phi||_{[L^2(\Omega)]^n}\bigr\}\cr}\right.\leqno(3.9)$$ The constants
do not depend on $\varepsilon$.\fin
\noindent{\bf Theorem 3.4 : }{\it For any  $\phi\in H^1(\Omega)$, there exists
$\widehat{\phi}_\varepsilon\in H^1_{per}(Y ; L^2(\Omega))$ such that
$$\left\{\eqalign{&||\widehat{\phi}_\varepsilon||_{H^1(Y ; L^2(\Omega))}\le C||\nabla_x\phi||_{[L^2(\Omega)]^n},\cr &||{\cal
T}_\varepsilon(\nabla_x\phi)-\nabla_x\phi-\nabla_y\widehat{\phi}_\varepsilon||_{ [L^2(Y ; H^{-1}(\Omega))]^n}\le
C\varepsilon||\nabla_x\phi||_{[L^2(\Omega)]^n}.}\right.\leqno(3.10)$$ 
\noindent  The constants  depend only on $n$ and 
$\partial\Omega$.}

\noindent{\bf Proof : } Let $\phi\in H^1(\Omega)$.  The function $\phi$ is decomposed 
$$\phi=\Phi+\varepsilon \underline{\phi} , \quad\hbox{where}\enskip \Phi={\cal Q}_\varepsilon(\phi)\quad\hbox{and} \quad
\underline{\phi}={1\over \varepsilon}{\cal R}_\varepsilon(\phi),\qquad {\cal R}_\varepsilon(\phi)=\phi-{\cal
Q}_\varepsilon(\phi),$$ with the following estimate (see [4]):
$$||\nabla_x\Phi||_{[L^2(\Omega)]^n}+||\underline{\phi}||_{L^2(\Omega)}+\varepsilon||\nabla_x
\underline{\phi}||_{[L^2(\Omega)]^n}\le C||\nabla_x\phi||_{[L^2(\Omega)]^n}.\leqno(3.11)$$ Proposition  3.3 applied to
$\underline{\phi}$ gives us the existence of an element
$\widehat{\phi}_{\varepsilon}$ in $H^1_{per}(Y ; L^2(\Omega))$ such that
$$\left\{\eqalign{  & ||\widehat{\phi}_{\varepsilon}||_{H^1(Y ; L^2(\Omega))}\le C||\nabla_x\phi||_{[L^2(\Omega)]^n},\cr
 & ||{\cal T}_\varepsilon(\underline{\phi})-\widehat{\phi}_{\varepsilon}||_{H^1(Y ;  H^{-1}(\Omega))}\le
C\varepsilon||\nabla_x\phi||_{[L^2(\Omega)]^n}.\cr}\right.\leqno(3.12)$$

\noindent   We  evaluate $||{\cal T}_\varepsilon(\nabla_x\Phi)-\nabla_x\Phi||_{[L^2(Y ;  H^{-1}(\Omega))]^n}$. 
\vskip 1mm
\noindent From the   inequality   $(3.2)$, applied to each partial derivative of $\Phi$, it follows
$$\Bigl\|{\partial\Phi\over \partial x_i}-M^\varepsilon_Y\Bigl({\partial\Phi\over \partial x_i}\Bigr)\Bigr\|_{ H^{-1}(\Omega)}\le
C\varepsilon||\nabla_x\Phi||_{[L^2(\Omega)]^n}\le C\varepsilon||\nabla_x\phi||_{[L^2(\Omega)]^n}\leqno(3.13)$$ There
results, from the definition of $\Phi$, that $\displaystyle y\longrightarrow{\cal T}_\varepsilon\Bigl({\partial\Phi\over
\partial x_i}\Bigr)(.,y)$ is linear with respect to each variable.  For any $\psi\in H^1_0(\Omega)$, we have
$$\eqalign{
<{\cal T}_\varepsilon \Bigl({\partial\Phi\over \partial x_1}\Bigr)(., y)-M^\varepsilon_Y\Bigl({\partial\Phi\over\partial
x_1}\Bigr),\psi>_{H^{-1}(\Omega)\, ,\, H^1_0(\Omega)}&=
\int_\Omega\Bigl\{{\cal T}_\varepsilon\Bigl({\partial\Phi\over \partial x_1}\Bigr)(., y)-M^\varepsilon_Y\Bigl({\partial\Phi\over
\partial x_1}\Bigr)\Bigr\}\psi\cr &=\int_{\Omega_\varepsilon}\Bigl\{{\cal T}_\varepsilon\Bigl({\partial\Phi\over \partial
x_1}\Bigr)(., y)-M^\varepsilon_Y\Bigl({\partial\Phi\over
\partial x_1}\Bigr)\Bigr\}M^\varepsilon_Y(\psi)\cr}$$ Set for 
$i=(i_1,\ldots, i_n)\in \{0,1\}^n$,$$ \overline{y}_i^{i_k}=\left\{\eqalign{ &y_k\qquad\hbox{if}\enskip i_k=1 \cr &1-
y_k\enskip \hbox{if}\enskip i_k=0 \cr}\right.$$ We have
$${\cal T}_\varepsilon\bigl(\Phi\bigr)(\varepsilon\xi, y)=\sum_{i_1\, \ldots, i_n}M^\varepsilon_Y(\phi)\bigl(
\varepsilon\xi+\varepsilon i\bigr)\overline{y}_1^{i_1}\ldots \overline{y}_n^{i_n},\qquad \xi=\Bigl[{x\over \varepsilon}\Bigr]$$
hence
$$\eqalign{ {\cal T}_\varepsilon\Bigl({\partial \Phi\over \partial x_1}\Bigr)(\varepsilon\xi, y)&=\sum_{i_2\, \ldots, i_n}{
M^\varepsilon_Y(\phi)\bigl(\varepsilon\xi+\varepsilon(1,i_2,...,i_n)\bigr)-M^\varepsilon_Y(\phi)\bigl(
\varepsilon\xi+\varepsilon(0,i_2,...,i_n)\bigr)\over \varepsilon}\overline{y}_2^{i_2}... \overline{y}_n^{i_n}\cr
M^\varepsilon_Y\Bigl({\partial \Phi\over \partial x_1}\Bigr)(\varepsilon\xi)&={1\over 2^{n-1}}\sum_{i_2\, \ldots,
i_n}{M^\varepsilon_Y(
\phi)\bigl(\varepsilon\xi+\varepsilon(1,i_2,\ldots,i_n)\bigr)-M^\varepsilon_Y(\phi)\bigl(
\varepsilon\xi+\varepsilon(0,i_2,\ldots,i_n)\bigr)\over \varepsilon} \cr}$$ We deduce that
$$\eqalign{
&\int_{\Omega_\varepsilon}\Bigl\{{\cal T}_\varepsilon\Bigl({\partial \Phi\over \partial x_1}\Bigr)(. , y)-
M^\varepsilon_Y\Bigl({\partial \Phi\over \partial x_1}\Bigr) \Bigr\}M^\varepsilon_Y(\Psi)=\cr
\varepsilon^n\sum_\xi\sum_{i_2\,..., i_n}&
\Bigl({M^\varepsilon_Y(\phi)\bigl(\varepsilon\xi+\varepsilon(1,i_2,...,i_n)\bigr)-M^\varepsilon_Y(\phi)\bigl(
\varepsilon\xi+\varepsilon(0,i_2,...,i_n)\bigr)\over \varepsilon} - M^\varepsilon_Y\Bigl({\partial \Phi\over \partial
x_1}\Bigr)(\varepsilon\xi)\Bigr)\cr &\times \overline{y}_2^{i_2}\ldots \overline{y}_n^{i_n} 
M^\varepsilon_Y(\psi)(\varepsilon\xi)\cr}$$ The above integral is equal to
$$\varepsilon^n\sum_\xi {M^\varepsilon_Y(\phi)\bigl(\varepsilon\xi+\varepsilon\vec e_1\bigr)-M^\varepsilon_Y(\phi) 
(\varepsilon\xi)\over \varepsilon} \sum_{i_2\, \ldots, i_n}
\bigl(M^\varepsilon_Y(\psi)\bigl(\varepsilon\xi-\varepsilon(0,i_2,..., i_n)\bigr )-\overline{M^\varepsilon_Y}\bigl( \psi\bigr)
(\varepsilon\xi)\Bigr) \overline{y}_2^{i_2}...\overline{y}_n^{i_n} $$ where
$$\overline{M^\varepsilon_Y}\bigl( \psi\bigr)(\varepsilon\xi)={1\over 2^{n-1}}\sum_{i_2,\ldots,i_n}M^\varepsilon_Y(\psi) 
\bigl(\varepsilon\xi-\varepsilon(0,i_2,\ldots, i_n)\bigr )$$ which gives the following inequality
$$ <{\cal T}_\varepsilon \Bigl({\partial\Phi\over \partial x_1}\Bigr)(., y)-M^\varepsilon_Y\Bigl({\partial\Phi\over\partial
x_1}\Bigr),\psi>_{H^{-1}(\Omega)\, ,\, H^1_0(\Omega)}\le C\varepsilon\bigl|y_2^{i_2}\ldots
 y_n^{i_n}\bigr|||\nabla_x\phi||_{[L^2(\Omega)]^n}||\psi||_{H^1_0(\Omega)}$$and  
$$\forall y\in Y,\qquad \Bigl\|{\cal T}_\varepsilon\Bigl({\partial \Phi\over \partial x_1}\Bigr)(., y)- M^\varepsilon_Y\Bigl({\partial
\Phi\over \partial x_1}\Bigr) \Bigr\|_{H^{-1}(\Omega)}
\le C\varepsilon ||\nabla_x\phi||_{[L^2(\Omega)]^n}.$$ Considering $(3.13)$ and all the partial derivatives, we obtain
$$||{\cal T}_\varepsilon(\nabla_x\Phi)-\nabla_x\Phi||_{[L^2(Y ; H^{-1}(\Omega))]^n}\le
C\varepsilon||\nabla_x\phi||_{[L^2(\Omega)]^n}$$ Thanks to $(3.12)$,  and to the above inequality and, moreover, to
$$||\varepsilon\nabla_x\underline{\phi}||_{[H^{-1}(\Omega)]^n}\le C\varepsilon||\underline{\phi}||_{L^2(\Omega)}\le
C\varepsilon||\nabla_x\phi||_{[L^2(\Omega)]^n}$$ the second estimate of $(3.10)$ is proved.\fin
\medskip
\noindent{\bf 4. Error estimate}
\medskip
\noindent We consider the following homogenization problem: find $\phi^{\varepsilon}\in H^1_{\Gamma_0}(\Omega)$ such that
$$\left\{\eqalign{&\forall \psi\in H^1_{\Gamma_0}(\Omega)=\bigl\{\phi\in H^1(\Omega)\; |\;
\phi=0\; \hbox{on}\;
\Gamma_0\bigr\},\cr&\int_\Omega A\bigl(\bigl\{{.\over\varepsilon}\bigr\}\bigr)\nabla\phi^{\varepsilon}.\nabla\psi=
\int_{\Omega}f\psi,}\right.\leqno(4.1)$$  where $\Omega$ is a bounded domain in $\R^n$ with lipschitzian boundary,
$\Gamma_0$ is a part of $\partial \Omega$ whose measure is nonnull or empty,  $f$ belongs to $L^p(\Omega)$,
$p>\displaystyle{2n\over n+2}$, (if $\Gamma_0=\emptyset$, we suppose that $\displaystyle \int_\Omega f=0$) and $A$ is a
square matrix of elements belonging to $L^\infty_{per}(Y)$, verifying the condition of uniform ellipticity  $c|\xi|^2\le
A(y)\xi.\xi\le C|\xi|^2$  a.e. $y\in Y$, with $c$ and  $C$ strictly positive constants. 

We have shown, see [4], that
$\nabla_x\phi^\varepsilon-\nabla_x\Phi-{\cal U}_\varepsilon\bigl(\nabla_y
\widehat{\phi}\bigr)$  strongly converges towards $0$ in
$[L^2(\Omega)]^n$, where ${\cal U}_\varepsilon$ is the averaging operator  defined by
$$\Psi\in L^2(\Omega\times Y)\qquad {\cal U}_\varepsilon(\Psi)(x)=\int_Y\Psi\Bigl(\varepsilon\Bigr[{x\over
\varepsilon}\Bigr]+\varepsilon z, \Bigl\{{x\over
\varepsilon}\Bigr\}\Bigr]\Bigr)dz,\qquad {\cal U}_\varepsilon(\Psi)\in L^2(\Omega),$$ and where  
$$(\Phi,\widehat{\phi})\in  H^1_{\Gamma_0}(\Omega)\times L^{2}(\Omega,H^1_{per}(Y)/\R)$$ is the solution of the limit
problem of unfolding homogenization
 $$\left\{\eqalign{
 &\forall (\Psi,\widehat{\psi})\in H^1_{\Gamma_0}(\Omega)\times L^{2}(\Omega ; H^1_{per}(Y)/\R)\cr&
\int_{\Omega}\int_{Y}A\bigl\{\nabla_x\Phi+\nabla_y\widehat{\phi}\bigr\}.
\bigl\{\nabla_x\Psi+\nabla_y\widehat{\psi}
\bigr\}=\int_\Omega f\Psi.}\right.\leqno(4.2)$$ 
\vskip-0.2cm
\noindent If $\Gamma_0=\emptyset$, we take $\displaystyle \int_\Omega\phi^\varepsilon=\int_\Omega\Phi=0$.

\noindent   We recall that the correctors $\chi_i$,
$i\in\{1,\ldots,n\}$, are the solutions of the following variational problems
\vskip-0.2cm
$$\chi_i\in H^1_{per}(Y),\qquad \int_YA(y)\nabla_y\bigl(\chi_i(y)+y_i\bigr)\nabla_y\psi(y)dy=0,\qquad
\forall\psi\in H^1_{per}(Y)$$
\noindent They allow us to express $\widehat{\phi}$ in terms of $\nabla_x\Phi$ 
$$\widehat{\phi}=\sum_{i=1}^n{\partial\Phi\over\partial x_i}\chi_i.$$

In Theorem 3  of [6]  our hypothesis was that the solution $\Phi$ of the homogenized problem belonged to $W^{2,p}(\Omega)$
($p>n$)  and we gave the following error estimate :
$$||\phi^\varepsilon-\Phi||_{L^2(\Omega)}+||\nabla_x\phi^\varepsilon-\nabla_x\Phi- \sum_{i=1}^n
{\partial\Phi\over\partial x_i}\nabla_y\chi_i\Bigl(\Bigl\{{.\over
\varepsilon}\Bigr\}\Bigr)||_{[L^2(\Omega)]^n}\le C\varepsilon^{\inf\{1/2,1-{n/p}\}},$$ the constant
depends on $n$, $p$, $A$, $||\Phi||_{W^{2,p}(\Omega)}$ and $\partial\Omega$. 
Then in Theorem 4 from [6]  we  obtained, by an interpolation method, the error estimate in the case where
$\Gamma_0=\partial\Omega$, and where the boundary of $\Omega$ is of class ${\cal C}^{1,1}$  and where $f$ belongs to
$L^p(\Omega)$ ($p\displaystyle >{2n\over n+2}$). 
$$||\phi^\varepsilon-\Phi||_{L^2(\Omega)}+||\nabla_x\phi^\varepsilon-\nabla_x\Phi- {\cal U}_\varepsilon\bigl(
\nabla_y\widehat{\phi}\bigr)||_{[L^2(\Omega)]^n}\le C\varepsilon^{\inf\{1/2,{1/2+1/n-1/p\over
1+1/n}\}}||f||_{L^p(\Omega)},$$ the constant depends on $n$, $p$, $A$ and $\partial\Omega$.

 If $\Phi$ belongs to $H^2(\Omega)$ the function $\displaystyle \sum_{i=1}^n {\partial\Phi\over\partial
x_i}\chi_i\Bigl(\Bigl\{{.\over\varepsilon}\Bigr\}\Bigr)$ does not  generally belong to $H^1(\Omega)$. However if $\Phi$ belongs 
only to
$H^1(\Omega)$ then from its definition the function
$\displaystyle {\cal Q}_\varepsilon\Bigl({\partial\Phi\over\partial x_i}\Bigr)$ belongs to $W^{1,\infty}(\Omega)$. Hence 
function $\displaystyle  {\cal Q}_\varepsilon\Bigl({\partial\Phi\over\partial x_i}\Bigr)
\chi_i\Bigl(\Bigl\{{.\over \varepsilon}\Bigr\}\Bigr)$ belongs to $H^1(\Omega)$ and thanks to $(3.5)$ it verifies 
$$\eqalign{
||{\cal Q}_\varepsilon\Bigl({\partial\Phi\over\partial
x_i}\Bigr)\chi_i\Bigl(\Bigl\{{.\over\varepsilon}\Bigr\}\Bigr)||_{L^2(\Omega)}&\le C
||\nabla_x\Phi||_{[L^2(\Omega)]^n}||\chi_i||_{L^2(Y)}\le C ||\nabla_x\Phi||_{[L^2(\Omega)]^n}\cr ||{\cal
Q}_\varepsilon\Bigl({\partial\Phi\over\partial x_i}\Bigr)\chi_i\Bigl(\Bigl\{{.\over\varepsilon}\Bigr\}\Bigr)||_{H^1(\Omega)}&\le
{C\over
\varepsilon} ||\nabla_x\Phi||_{[H^1(\Omega)]^n}||\chi_i||_{H^1(Y)}\le {C\over \varepsilon}
||\nabla_x\Phi||_{[H^1(\Omega)]^n}\cr}$$ This is the  reason why in the approximate solution we replace  $\displaystyle 
{\partial\Phi\over\partial x_i}$ with $\displaystyle {\cal Q}_\varepsilon\Bigl( {\partial\Phi\over\partial x_i} \Bigr)$. In the
following theorems we are going to obtain estimates that are better than those obtained in [6], with weaker hypotheses.
\medskip
\noindent{\bf 4.1 First case : Homogeneous Dirichlet or Neumann condition and boundary of class ${\cal C}^{1,1}$.}
\medskip
\noindent{\bf Theorem 4.1 : }{\it We suppose that  $\Omega$ is a ${\cal C}^{1,1}$ bounded domain in $\R^n$, 
$\Gamma_0=\partial\Omega$ and $f\in L^2(\Omega)$. Then we have
$$||\phi^\varepsilon-\Phi||_{L^2(\Omega)}+||\nabla_x\phi^\varepsilon-\nabla_x\Phi- \sum_{i=1}^n
 {\cal Q}_\varepsilon\Bigl({\partial\Phi\over\partial
x_i}\Bigr)\nabla_y\chi_i\Bigl(\Bigl\{{.\over\varepsilon}\Bigr\}\Bigr)||_{[L^2(\Omega)]^n}\le
C\varepsilon^{1/2}||f||_{L^2(\Omega)}.\leqno(4.3)$$ The constant depends on $n$,
$A$ and $\partial\Omega$.}

\medskip
\noindent{\bf Theorem 4.2 : }{\it We suppose that  $\Omega$ is a ${\cal C}^{1,1}$ bounded domain in $\R^n$,
$\Gamma_0=\emptyset$, $f\in L^2(\Omega)$. Then we have
$$||\phi^\varepsilon-\Phi||_{L^2(\Omega)}+||\nabla_x\phi^\varepsilon-\nabla_x\Phi-\sum_{i=1}^n {\cal Q}_\varepsilon\Bigl(
{\partial\Phi\over\partial x_i}\Bigr)\nabla_y\chi_i\Bigl(\Bigl\{{.\over\varepsilon}\Bigr\}\Bigr)||_{[L^2(\Omega)]^n}\le C
\varepsilon^{1/ 2}||f||_{L^2(\Omega)}\leqno(4.4)$$ The constant depends on  $n$, $A$ and $\partial\Omega$.}
\medskip
\noindent The proof of Theorems 4.1 and 4.2 is based on the following proposition.
\medskip
\noindent{\bf Proposition 4.3 : }{\it We suppose that the solution $\Phi$ of the unfolded problem belongs to
$H^2(\Omega)$. Therefore we have
$$||\phi^\varepsilon-\Phi||_{L^2(\Omega)}+||\nabla_x\phi^\varepsilon-\nabla_x\Phi-\sum_{i=1}^n {\cal Q}_\varepsilon\Bigl(
{\partial\Phi\over\partial x_i}\Bigr)\nabla_y\chi_i\Bigl(\Bigl\{{.\over\varepsilon}\Bigr\}\Bigr)||_{[L^2(\Omega)]^n}\le C
\varepsilon^{1/ 2}\leqno(4.5)$$  The
constant depends on $A$, $n$, $||\Phi||_{H^2(\Omega)}$ and $\partial\Omega$.}

\noindent{\bf Proof : }   We denote by  $\rho(x)=dist(x,\partial \Omega)$ the distance  between $x\in \Omega$ and the boundary
of $\Omega$.

\noindent We show that if $(\Phi,\widehat{\phi})$ is the solution of the unfolded problem, then
 $\displaystyle \Phi+\sum_{i=1}^n\varepsilon\rho_\varepsilon{\cal Q}_\varepsilon\Bigl({\partial\Phi\over
\partial x_i}\Bigr)\chi_i\Bigl(\Bigl\{{.\over\varepsilon}\Bigr\}\Bigr)$ is an approximate solution to the homogenization
 problem $(4.1)$; $\displaystyle\rho_\varepsilon(.)=\inf\Bigl\{{\rho(.)\over\varepsilon},1\Bigr\}$. The presence of the function
$\rho_\varepsilon$ in the sum guarantees the nullity of the  approximate solution on $\Gamma_0$.
\vskip 1mm
\noindent{\bf Step one. } We present some estimates of $\rho_\varepsilon$, $\nabla_x\Phi$ and
$\displaystyle\chi_i\Bigl(\Bigl\{{.\over\varepsilon}\Bigr\}\Bigr)$ on  the neighborhood 
$\widehat{\Omega}_\varepsilon=\bigl\{x\in\Omega\; |\; \rho(x)<\varepsilon\bigr\}$ of the boundary of  $\Omega$. We have
$$\left\{\eqalign{
& ||\nabla_x\rho_\varepsilon||_{[L^\infty(\Omega)]^n}=||\nabla_x\rho_\varepsilon||_{[L^\infty(\widehat{
\Omega}_\varepsilon)]^n}= \varepsilon^{-1},\cr  
&||\nabla_x\Phi||_{[L^2(\widehat{\Omega}_\varepsilon)]^n}\le C\varepsilon^{1/2}||\Phi||_{H^2(\Omega)}\cr
&\Longrightarrow \;\, ||{\cal Q}_\varepsilon(\nabla_x\Phi)||_{[L^2(\widehat{
\Omega}_\varepsilon)]^n}+||M^\varepsilon_Y(\nabla_x\Phi)||_{[L^2(\widehat{\Omega}_\varepsilon)]^n}\le
C\varepsilon^{1/2}||\Phi||_{H^2(\Omega)},\cr
&\big\|\chi_i\Bigl(\Bigl\{{.\over\varepsilon}\Bigr\}\Bigr)\Bigr\|_{L^2(\widehat{\Omega}_\varepsilon)}+\Big\|\nabla_y\chi_i\Bigl(
\Bigl\{{.\over\varepsilon}\Bigr\}\Bigr)\Bigr\|_{[L^2(\widehat{\Omega}_\varepsilon)]^n}\le
C\varepsilon^{ 1/2}||\nabla_y\chi_i||_{[L^2(Y)]^n}\le C\varepsilon^{ 1/2}.\cr}\right.\leqno(4.6)$$ 
\noindent The estimate  of  $\rho_\varepsilon$ follows from its definition.  The estimate  of
$\nabla_x\Phi$ in $[L^2(\widehat{\Omega}_\varepsilon)]^n$  comes from the gradient belonging to $H^2(\Omega)$.
 The number of  cells covering $\widehat{\Omega}_\varepsilon$ is of  order of
$\varepsilon^{1-n}$, hence we obtain the estimates of $\nabla_y\chi_i$ and $\chi_i$ on the neighborhood of the boundary of
$\Omega$. We will note for the rest of the demonstration that the
support  of $1-\rho_\varepsilon$ is contained in $\Omega\setminus\widehat{\Omega}_\varepsilon$.
\vskip 1mm
\noindent{\bf Step two. } Let $\Psi\in H^1_{\Gamma_0}(\Omega)$. Thanks to Theorem 3.4, there exists
$\widehat{\psi}_\varepsilon\in H^1_{per}(Y ;  L^2(\Omega))$ verifying the estimates $(3.10)$. We take the couple
$(\Psi,\widehat{\psi}_\varepsilon)$ as a test-function in the unfolded problem  $(4.2)$ and we introduce
$\rho_\varepsilon$. The gradient of $\Phi$ belongs to $[H^1(\Omega)]^n$, and according to $(4.6)$
$$||(1-\rho_\varepsilon)\nabla_x\Phi||_{[L^2(\Omega)]^n}\le ||\nabla_x\Phi||_{[L^2(\widehat{\Omega}_\varepsilon)]^n}
\le C\varepsilon^{1/ 2}||\Phi||_{H^2(\Omega)},\leqno(4.7)$$ which gives us
$$\Bigl|\int_\Omega f\Psi-\int_{\Omega\times Y}A(y)\rho_\varepsilon(x)\Bigl\{\nabla_x\Phi(x)+\sum_{i=1}^n{\partial\Phi\over
\partial x_i}(x)\nabla_y\chi_i(y)\Bigr\}\bigl(\nabla_x\Psi+\nabla_y\widehat{\psi}_\varepsilon\bigr)\Bigr|\le C\varepsilon^{1/
2}||\Psi||_{H^1(\Omega)}.$$ 
 \noindent In the integral on $\Omega\times Y$ we replace $\nabla_x\Psi+\nabla_y\widehat{\psi}_\varepsilon$ by ${\cal
T}_\varepsilon(\nabla_x\Psi)$,  thanks to  $(3.10)$ of  Theorem 3.4.  The function $\rho_\varepsilon\nabla_x\Phi$ belongs to
$[H^1_0(\Omega)]^n$ and verifies
$||\rho_\varepsilon\nabla_x\Phi||_{ [H^1(\Omega)]^n}\le C\varepsilon^{-{1/ 2}}||\Phi||_{H^2(\Omega)}$ for
$$\eqalign{ &\Bigl\|\nabla_x\Bigl\{\rho_\varepsilon{\partial\Phi\over
\partial x_i}\Bigr\}\Bigr\|_{[L^2(\Omega)]^n}\le
\Bigl\|\nabla_x\rho_\varepsilon {\partial\Phi\over \partial x_i}\Bigr\|_{[L^2(\Omega)]^n}+\Bigl\|\rho_\varepsilon
\nabla_x\Bigl\{{\partial\Phi\over \partial x_i}\Bigr\}\Bigr\|_{[L^2(\Omega)]^n}\cr
\le & ||\nabla_x\rho_\varepsilon||_{[L^\infty(\widehat{\Omega}_\varepsilon)]^n}\Bigl\|{\partial\Phi\over
\partial x_i}\Bigr\|_{L^2(\widehat{\Omega}_\varepsilon)}+\Bigl\|\nabla_x\Bigl\{{\partial\Phi\over
\partial x_i}\Bigr\}\Bigr\|_{[L^2(\Omega)]^n}\le C\varepsilon^{-{1/ 2}}||\Phi||_{H^2(\Omega)}.\cr}$$ Then we remove
$\rho_\varepsilon$ in the products
$\rho_\varepsilon(x)\nabla_x\Phi(x)$ and $\displaystyle
\rho_\varepsilon(x){\partial\Phi\over\partial x_i}(x)\nabla_y\chi_i(y)$ by using $(4.7)$ again. And then we replace
 $\nabla_x\Phi$ with $M^\varepsilon_Y(\nabla_x\Phi)$ and  in the sum  we replace $\displaystyle {\partial\Phi\over\partial x_i}$
with $\displaystyle M^\varepsilon_Y\Bigl({\partial\Phi\over\partial x_i}\Bigr)$.  Thanks to $(3.2)$, we obtain
$$\eqalign{  &\Bigl|\int_\Omega f\Psi-{1\over |Y|}\int_{\Omega\times Y}A(.)\Bigl\{M^\varepsilon_Y(\nabla_x\Phi)+
\sum_{i=1}^nM^\varepsilon_Y\Bigl({\partial\Phi\over\partial x_i}\Bigr)\nabla_y\chi_i(.)\Bigr\}{\cal
T}_\varepsilon(\nabla_x\Psi)\Bigr|\le C\varepsilon^{1/2}||\Psi||_{H^1(\Omega)}\cr}$$ By inverse unfolding we transform  the
integral on
$\Omega\times Y$ into an  integral on $\Omega$. Then  we replace  $M^\varepsilon_Y(
\nabla_x\Phi)$ with $\nabla_x\Phi$  and we  reintroduce
$\rho_\varepsilon$ in front of $\displaystyle M^\varepsilon_Y\Bigl({\partial\Phi\over\partial
x_i}\Bigr)\nabla_y\chi_i\Bigl(\Bigl\{{.\over\varepsilon}\Bigr\}\Bigr)$ 
$$\Bigl\|(1-\rho_\varepsilon)M^\varepsilon_Y\Bigl({\partial\Phi\over\partial
x_i}\Bigr)\nabla_y\chi_i\Bigl(\Bigl\{{.\over\varepsilon}\Bigr\}\Bigr)\Bigr\|_{L^2(\Omega)} \le
\Bigl\|M^\varepsilon_Y\Bigl({\partial\Phi\over\partial x_i}\Bigr)\Bigr\|_{L^2(\widehat{\Omega
}_\varepsilon)}\Bigl\|\nabla_y\chi_i\Bigl(\Bigl\{{.\over\varepsilon}\Bigr\}\Bigr)
\Bigr\|_{[L^2(\widehat{\Omega }_\varepsilon)]^n} \le C\varepsilon ||\Phi||_{H^2(\Omega)}.\leqno(4.8)$$ This done, we have
$$\Bigl|\int_\Omega f\Psi-\int_{\Omega}A\Bigl({.\over\varepsilon}\Bigr)\Bigl\{\nabla_x\Phi+
\sum_{i=1}^n\rho_\varepsilon M^\varepsilon_Y\Bigl({\partial\Phi\over\partial
x_i}\Bigr)\nabla_y\chi_i\Bigl(\Bigl\{{.\over\varepsilon}\Bigr\}\Bigr)\Bigr\}\nabla_x\Psi\Bigr| \le 
C\varepsilon^{1/2}||\Psi||_{H^1(\Omega)}. $$  From $(3.4)$  we obtain
$$\Bigl\|\Bigl\{M^\varepsilon_Y\Bigl({\partial\Phi\over\partial x_i}\Bigr)-{\cal Q}_\varepsilon\Bigl({\partial\Phi\over\partial
 x_i}\Bigr)\Bigr\}\nabla_y\chi_i\Bigl(\Bigl\{{.\over\varepsilon}\Bigr\}\Bigr)\Bigr\|_{L^2(\Omega)}\le C\varepsilon
||\Phi||_{H^2(\Omega)}$$ hence
$$\Bigl|\int_\Omega
f\Psi-\int_{\Omega}A\Bigl(\Bigl\{{.\over\varepsilon}\Bigr\}\Bigr)\Bigl\{\nabla_x\Phi+\sum_{i=1}^n\rho_\varepsilon
 {\cal Q}_\varepsilon\Bigl({\partial\Phi\over\partial
x_i}\Bigr)\nabla_y\chi_i\Bigl(\Bigl\{{.\over\varepsilon}\Bigr\}\Bigr)\Bigr\}\nabla_x\Psi\Bigr| \le 
C\varepsilon^{1/2}||\Psi||_{H^1(\Omega)}. $$
\noindent We now estimate the terms which appear in the calculation of the gradient of the approximate solution but do not
appear in the above expression. Thanks to  $(4.6)$ and $(3.5)$ we have
$$\left\{\eqalign{
\Big\|\varepsilon{\partial\rho_\varepsilon\over \partial x_j} {\cal Q}_\varepsilon\Bigl({\partial\Phi\over\partial x_i}\Bigr)
\chi_i\Bigl(\Bigl\{{.\over\varepsilon}\Bigr\}\Bigr)\Bigr\|_{L^2(\Omega)}  &\le
\Big\|\varepsilon{\partial\rho_\varepsilon\over \partial x_j}\Bigr\|_{L^\infty(\widehat{\Omega}_\varepsilon)}
\Big\| {\cal Q}_\varepsilon\Bigl({\partial\Phi\over\partial x_i}\Bigr)\Bigr\|_{L^2(\widehat{\Omega}_\varepsilon)}
||\chi_i||_{L^2(Y)}\cr&\vphantom{1\over 2}\le C\varepsilon^{1/2} ||\Phi||_{H^2(\Omega)},\cr
\Big\|\varepsilon\rho_\varepsilon{\partial  \over\partial x_j} {\cal Q}_\varepsilon\Bigl({\partial\Phi\over\partial
x_i}\Bigr)\chi_i\Bigl(\Bigl\{{.\over\varepsilon}\Bigr\}\Bigr)
\Bigr\|_{L^2(\Omega)}  &\le \varepsilon||\rho_\varepsilon||_{L^\infty(\Omega)}\Big\|\nabla {\cal Q}_\varepsilon
\Bigl({\partial\Phi\over\partial x_i} \Bigr) \Bigr\|_{L^2(\Omega)}||\chi_i||_{L^2(Y)}\cr&\vphantom{1\over 2} \le 
C\varepsilon||\Phi||_{H^2(\Omega)}.\cr}\right.\leqno(4.9)$$  
\noindent Now we use the equality $$\displaystyle
\int_\Omega f\Psi=\int_{\Omega}A\Bigl(\Bigl\{{x\over
\varepsilon}\Bigr\}\Bigr)\nabla_x\phi^\varepsilon(x)\nabla_x\Psi(x),$$ and we take as a test function
$$\Psi=\displaystyle\phi^\varepsilon-\Bigl[\Phi+\sum_{i=1}^n\varepsilon\rho_{\varepsilon}{\cal Q}_\varepsilon\Bigl(
{\partial\Phi\over \partial x_i}\Bigr)\chi_i\Bigl(\Bigl\{{.\over\varepsilon}\Bigr\}\Bigr)\Bigr],$$ to obtain
$$\Bigl\|\nabla_x\phi^\varepsilon-\nabla_x\Bigl[\Phi+\sum_{i=1}^n\varepsilon\rho_\varepsilon{\cal
Q}_\varepsilon\Bigl({\partial
\Phi\over\partial x_i}\Bigr)\chi_i\Bigl(\Bigl\{{.\over\varepsilon}\Bigr\}\Bigr)\Bigr]\Bigr\|_{[L^2(\Omega)]^n} \le C
\varepsilon^{1/2}$$ This gives the estimate  $(4.5)$ thanks  to the Poincar\'e inequality or the Poincar\'e-Wirtinger 
inequality and $(4.7)$  and $(4.9)$.\fin
\medskip
\noindent{\bf  Proofs of Theorems 4.1 and 4.2  : } The boundary of $\Omega$ is of class   ${\cal C}^{1,1}$, then for  $f\in
L^2(\Omega)$,  the solution  $\Phi$ from
$(4.2)$ with
$\Gamma_0=\partial\Omega$ or $\Gamma_0=\emptyset$, belongs to
$H^2(\Omega)$  and verifies $\displaystyle ||\Phi||_{H^2(\Omega)}\le C||f||_{L^2(\Omega)}$. \fin
\medskip
\noindent{\bf Corollary of Theorem 4.1 : }  When the correctors belong to
$W^{1,\infty}(Y)$ we obtain the classical error estimate (see $[2]$, $[5]$ and $[8]$).\fin
\medskip
\noindent{\bf 4.2 Second case : lipschitzian boundary.}
\medskip
\medskip\noindent{\bf Proposition 4.4 : }{\it We suppose that  solution $\Phi$ of the unfolded problem $(4.2)$ belongs to
$H^2_{loc}(\Omega)\cap W^{1,q}(\Omega)$ (
$ q>2$)  and verifies $$||\rho\nabla_x\Phi ||_{[H^1(\Omega)]^n}<+\infty,\leqno(4.10)$$
\noindent Then we have
$$||\phi^\varepsilon-\Phi||_{L^2(\Omega)}+ ||\nabla_x\phi^\varepsilon-\nabla_x\Phi-\sum_{i=1}^n {\cal
Q}_\varepsilon\Bigl({\partial\Phi\over
\partial x_i}\Bigr)\nabla_y\chi_i\Bigl(\Bigl\{{.\over\varepsilon}\Bigr\}\Bigr)||_{[L^2(\Omega)]^n}\le C\varepsilon^{{q-2\over
3q-2}}
\leqno(4.11)$$ The constant depends on $A$,
$n$,  $q$, $||\Phi||_{W^{1,q}(\Omega)}+||\rho\nabla_x\Phi||_{[H^1(\Omega)]^n}$  and $\partial\Omega$.}

\noindent{\bf Proof : }  We equip $W^{1,q}(\Omega)\cap H^2_{loc}(\Omega)$ with the norm
$$|||\Psi|||=||\Psi||_{W^{1,q}(\Omega)}+||\rho\nabla_x\Psi||_{[H^1(\Omega)]^n}$$
\noindent As in proposition 7, we show  that if $(\Phi,\widehat{\phi})$ is the solution of $(4.2)$, then
 $$\displaystyle\Phi+\sum_{i=1}^n\varepsilon\rho_{\varepsilon,\alpha}{\cal Q}_\varepsilon\Bigl({\partial\Phi\over
\partial x_i}\Bigr)\chi_i\Bigl(\Bigl\{{.\over\varepsilon}\Bigr\}\Bigr)$$ is an approximate solution of problem $(4.1)$, where 
$\displaystyle
\rho_{\varepsilon,\alpha}(.)=\inf\Bigl\{{d(.)\over\varepsilon^\alpha},1\Bigr\}$,
$\alpha$ belongs to interval
$]0,1]$ and will be fixed later.
\vskip 1mm
\noindent{\bf Step one. } We present some estimates of $\rho_{\varepsilon,\alpha}$  and $\Phi$ on  the
neighborhood $\widehat{\Omega}_{\varepsilon,\alpha}=\bigl\{x\in \Omega\; ; \; \rho(x)<\varepsilon^\alpha\bigr\}$ of the
boundary of $\Omega$. We have
$$\left\{\eqalign{ 
&{\vphantom{1\over 2}}||\nabla_x\rho_{\varepsilon,\alpha}||_{[L^\infty(\Omega)]^n}=
||\nabla_x\rho_{\varepsilon,\alpha}||_{[L^\infty(\widehat{\Omega}_{\varepsilon,\alpha})]^n}=\varepsilon^{-\alpha},\cr
&{\vphantom{1\over 2}}||\nabla_x\Phi||_{[L^2(\widehat{\Omega}_{\varepsilon,\alpha})]^n}\le C\varepsilon^{\alpha({1\over
2}-{1\over q})} ||\nabla_x\Phi||_{[L^q(\Omega)]^n},\cr &{\vphantom{1\over 2}}
||\rho_{\varepsilon,\alpha}\nabla_x\Phi||_{[H^1(\Omega)]^n}\le C\varepsilon^{-\alpha}|||\Phi|||.\cr}\right.\leqno(4.12)$$
\noindent{\bf Step two. } Let $\Psi\in H^1_{\Gamma_0}(\Omega)$. Thanks to Theorem 3.4, there exists
$\widehat{\psi}_\varepsilon\in L^2(\Omega ; H^1_{per}(Y))$ verifying the estimates $(3.10)$. We take the couple
$(\Psi,\widehat{\psi}_\varepsilon)$ as test-function in the unfolded  problem  $(4.2)$ and we introduce 
$\rho_{\varepsilon,\alpha}$.  The gradient of $\Phi$ verifies
$$ ||(1-\rho_{\varepsilon,\alpha})\nabla_x\Phi||_{[L^2(\Omega)]^n}\le
||\nabla_x\Phi||_{[L^2(\widehat{\Omega}_{\varepsilon,\alpha})]^n}\le C\varepsilon^{\alpha\{{1\over 2}-{1\over
q}\}}|||\Phi|||,\leqno(4.13)$$  according to $(4.12)$. This gives us
$$\Bigl|\int_\Omega f\Psi-\int_{\Omega\times Y}A\rho_{\varepsilon,\alpha}\Bigl\{\nabla_x\Phi+
\sum_{i=1}^n{\partial\Phi\over\partial x_i}\nabla_y\chi_i\Bigr\}\bigl(\nabla_x\Psi+\nabla_y
\widehat{\psi}_\varepsilon\bigr)\Bigr|\le C\varepsilon^{\alpha\{{1\over 2}-{1\over q}\}}||\Psi||_{H^1(\Omega)}$$ 
\noindent In the integral on $\Omega\times Y$  we replace  $\nabla_x\Psi+\nabla_y\widehat{
\psi}_\varepsilon$ with ${\cal T}_\varepsilon(\nabla_x\Psi)$,  thanks to $(3.10)$ from Theorem 3.4 and to $(4.12)$. Function
$\rho_{\varepsilon,\alpha}\nabla_x\Phi$ belongs to $[H^1_0(\Omega)]^n$ and thanks to $(3.2)$,  $(4.12)$ and $(4.13)$ we
get
$$\Bigl\|\rho_{\varepsilon,\alpha}{\partial\Phi\over \partial x_i}-\rho_{\varepsilon,\alpha}M^\varepsilon_Y\Bigl(
{\partial\Phi\over \partial x_i}\Bigr)\Bigr\|_{L^2(\Omega)} \le C\varepsilon^{\inf\{\alpha({1\over 2}-{1\over
q}),1-\alpha\}}|||\Phi||| $$ We also have $||{\cal T}_\varepsilon( \rho_{\alpha,\varepsilon})-\rho_{\varepsilon,\alpha}||_{
L^\infty(\Omega\times Y)}\le C\varepsilon^{1-\alpha}$. Now we proceed as in Proposition 4.3  to obtain
$$\Bigl|\int_\Omega f\Psi-\int_{\Omega}A\Bigl(\Bigl\{{.\over\varepsilon}\Bigr\}\Bigr)\Bigl\{\nabla_x\Phi+\sum_{i=1}^n\rho_{
\varepsilon,\alpha}{\cal Q}_\varepsilon\Bigl({\partial\Phi\over\partial
x_i}\Bigr)\nabla_y\chi_i\Bigl(\Bigl\{{.\over\varepsilon}\Bigr\}\Bigr)\Bigr\}\nabla_x\Psi\Bigr| 
\le  C\varepsilon^{\inf\{\alpha({1\over 2}-{1\over q}),1-\alpha\}}||\Psi||_{H^1(\Omega)}.$$ We choose
$\alpha=\displaystyle{2q\over 3q-2}$.
\noindent We  estimate the terms that appear in the calculation of the gradient of the approximate solution but do not
appear in the above expression thanks to  $(4.12)$. We now   use the equality $$\displaystyle \int_\Omega
f\Psi=\int_{\Omega}A\Bigl(\Bigl\{{.\over\varepsilon}\Bigr\}\Bigr)\nabla_x\phi^\varepsilon(x)\nabla_x\Psi(x)$$ and we take 
$$\Psi=\displaystyle\phi^\varepsilon-\Bigl(\Phi+\sum_{i=1}^n\varepsilon \rho_{\varepsilon,\alpha}{\cal
Q}_\varepsilon\Bigl({\partial\Phi\over \partial x_i}\Bigr)\chi_i\Bigl(\Bigl\{{.\over\varepsilon}\Bigr\}\Bigr)$$ as test-function,
to obtain
$$||\nabla_x\phi^\varepsilon-\nabla_x\Bigl(\Phi+\sum_{i=1}^n\varepsilon\rho_{\varepsilon,\alpha}{\cal Q}_\varepsilon\Bigl(
{\partial\Phi\over\partial x_i}\Bigr)\chi_i\Bigl(\Bigl\{{.\over\varepsilon}\Bigr\}\Bigr)\Bigr\|_{[L^2(\Omega)]^n} \le
C\varepsilon^{q-2\over 3q-2}.$$ \fin
\noindent{\bf Theorem 4.5 : }{\it We suppose that  $\Omega$ is a bounded domain in $\R^n$ with  lipschitzian boundary and
$\Gamma_0$ is a  union of connected components of $\partial\Omega$. Then, there exists $\gamma$ in the interval
$\displaystyle\Bigl]0,{1\over 3}\Bigr]$ depending on $A$, $n$ and $\partial \Omega$ such that for any $f\in L^2(\Omega)$
$$||\phi^\varepsilon-\Phi||_{L^2(\Omega)}+||\nabla_x\phi^\varepsilon-\nabla_x\Phi-\sum_{i=1}^n {\cal Q}_\varepsilon\Bigl(
{\partial\Phi\over\partial x_i}\Bigr)\nabla_y\chi_i\Bigl(\Bigl\{{.\over\varepsilon}\Bigr\}\Bigr)||_{[L^2(\Omega)]^n}
\le C\varepsilon^\gamma ||f||_{L^2(\Omega)}\leqno(4.14)$$ The constant
$C$ depends on $n$, $A$ and $\partial\Omega$.}

\noindent{\bf Proof : }  

\noindent {\bf  Step one. } We denote ${\cal A}$ the square matrix associated to the homogenized operator  (see
[5]). Let $R>0$ such that $\Omega\subset B(O; R)$ and $w\in H^1_0(B(O; R))$ the solution of the variational problem
$$\int_{B(O;R)}{\cal A}\nabla_x w\nabla_x v=\int_{\Omega}f v\qquad \forall v\in H^1_0(B(O;R))$$ 
We have $w\in H^2(B(O; R))$ and $||w||_{H^2(B(O;R)}\le C||f||_{L^2(\Omega)}$. The function $\Phi$ is solution
of the homogenized problem (see [5])
$$\int_{\Omega}{\cal A}\nabla_x \Phi\nabla_x v=\int_{\Omega}f v\qquad \forall v\in H^1_{\Gamma_0}(\Omega)$$ Hence $-div\bigl({\cal A}(\nabla_x
w-\nabla_x\Phi)\bigr)=0$ in $H^{-1}(\Omega)$ and $w-\Phi$ belongs to ${\cal C}^\infty(\Omega)$. We also have
$$\forall i\in\{1,\ldots, n\}\qquad\qquad -div\bigl({\cal A}(\nabla_x {\partial w\over \partial x_i}-\nabla_x{\partial \Phi\over \partial x_i})\bigr)=0\qquad
\hbox{in}\qquad {\cal D}^\prime(\Omega)$$   Lemma 2.2 of [7] gives us
$$ \forall i\in\{1,\ldots, n\}\qquad\qquad \int_\Omega \rho^2|\nabla_x({\partial w\over \partial x_i}-{\partial \Phi\over \partial x_i}) |^2  \le C\int_\Omega|{\partial
w\over \partial x_i}-{\partial \Phi\over \partial x_i} |^2$$  From the estimates of $w$ and $\Phi-w$, it follows:
$||\Phi||_{H^1(\Omega)}+||\rho\nabla_x\Phi||_{H^1(\Omega)}\le C||f||_{L^2(\Omega)}$.

\noindent{\bf Step two. }   Theorem A.3  of [3] asserts the existence of a real  $q>2$, depending on  $A$ and $\partial
\Omega$,
 such that $\Phi$ belongs to $W^{1,q}(\Omega)$. Thanks to   Proposition 4.4 we obtain  Theorem 4.5.\fin
\smallskip
\noindent{\bf Comments : }   In Theorems 4.1 and 4.2, if $\nabla\Phi=0$ on $\partial \Omega$, the error estimate is of order
$\varepsilon$.

\noindent In Theorem 4.1  if we replace $f\in L^2(\Omega)$ by  $f\in L^p(\Omega)$  $\Bigl(\displaystyle
{2n\over n+2}<p\le 2\Bigr)$, then  we prove that
$${\vphantom{1\over 2}}||\phi^\varepsilon-\Phi||_{L^2(\Omega)}+||\nabla_x\phi^\varepsilon-\nabla_x\Phi- \sum_{i=1}^n {\cal
Q}_\varepsilon\Bigl({\partial\Phi\over\partial
x_i}\Bigr)\nabla_y\chi_i\Bigl(\Bigl\{{.\over\varepsilon}\Bigr\}\Bigr)||_{[L^2(\Omega)]^n}\le C\varepsilon^{\inf\{{1\over 2},
n({1\over n}+{1\over 2}-{1\over p})\}}||f||_{L^p(\Omega)}.$$  The constant depends on  $n$,
$p$, $A$ and $\partial\Omega$. 

\noindent The estimates obtained in Theorems 4.1, 4.2 and 4.5 remain true if we suppose that the coefficients
of the square matrix $A$ of problem $(4.1)$  belong to $W^{1,\infty}(\Omega; L^\infty_{per}(Y))$.\fin
\vskip -1mm
\medskip\noindent{\bf Acknowledgements}
\smallskip
\noindent The author wishes to thank Do\"\i na Cioranescu and Alain Damlamian for their fruitful discussions.
 
\noindent The idea of the measurement of the periodic defect is due to Alain Damlamian.
\bigskip
\centerline{\bf References}
\medskip
\noindent [1] R. A. Adams Sobolev Spaces, Academic Press, New York, 1975.

\noindent [2] A. Bensoussan, J.-L.Lions and G.Papanicolaou, Asymptotic Analysis for Periodic Structures, North Holland,
Amsterdam, 1978.

\noindent [3]   M. Briane,  A. Damlamian and P. Donato,  
H-convergence in perforated domains, Nonlinear Partial Differential  
Equations and Their Applications, Coll\`ege de France seminar vol. XIII,
 D. Cioranescu and J.-L. Lions eds.,
 Pitman Research Notes in Mathematics Series 391, Longman, New York (1998), 
62--100.

\noindent [4] D. Cioranescu, A. Damlamian and  G. Griso, Periodic unfolding and homogenization, C. R. Acad. Sci. Paris, Ser. I 335
(2002), 99--104.

\noindent [5] D. Cioranescu and P. Donato,  An  Introduction to Homogenization. Oxford Lecture Series in Mathematics ans its
Applications 17, Oxford University Press, 1999.

\noindent [6]  G. Griso, Estimation d'erreur et \'eclatement en homog\'en\'eisation p\'eriodique.  C. R. Acad. Sci. Paris, Ser. I 335
(2002), 333--336.

\noindent [7]  G. Griso and B. Miara, Modelling of periodic electromagnetic  structures. Bianisotropic materials with  memory. (To appear)

\noindent [8]  O.A. Oleinik, A. S. Shamaev and G. A. Yosifian, Mathematical Problems in Elasticity and Homogenization,
North-Holland, Amsterdam, 1994. 

\bye